# Quantified 'shock-sensitivity' above the Maxwell load

J. M. T. Thompson[1], G. H. M. van der Heijden[2]

[1]*Department of Applied Maths & Theoretical Physics, University of Cambridge, CB3 0WA*
[2]*Department of Civil, Environmental and Geomatic Engineering, UCL, London, WC1E 6BT, UK*

**Abstract**

Using the static-dynamic analogy, work at Bath and Bristol has uncovered the vital organizing role of the Maxwell 'energy criterion' load in the advanced post-buckling of long-thin structures which exhibit severe shell-like imperfection sensitivity. It has become clear that above the Maxwell load, $P_M$, there are localized solutions offering an order-of-magnitude increase in sensitivity to lateral side-loads, whether static or dynamic. We propose to call this 'shock-sensitivity', and notice that so far only the seminal paper by Horak, Lord and Peletier in 2006 has quantified this in terms of an $E(P)$ energy-barrier versus load graph. In this paper we present three graphs of this nature for archetypal problems: the free twisted rod, the cylindrically constrained rod, and the strut on a softening elastic foundation. We find in all cases that the energy barrier of the localizing solution above $P_M$ is quite close to the energy of a *single* periodic wave. Now a single such wave is not kinematically admissible, and the corresponding periodic barrier must be for all the waves in the long structure, $N$, say. So in practice $N$ will be large, and does indeed tend to infinity with the length of the structure. Thus the sensitivity increases by a factor of a large $N$ as the Maxwell load is exceeded. This is important in its own right, and we do not seek to explain or fit curves to the scattered experimental buckling loads of shell structures.

**1. Introduction**

The buckling of thin aero-space shells has a long history. Experiments at the beginning of the 20[th] century showed that compressed spherical and cylindrical shells were failing at about one quarter of classical linear buckling loads, $P_C$. Von Karman and Tsien [1939, 1941] showed that there exists a very unstable, sub-critical post-buckling path of periodic equilibrium states, which falls rapidly from $P_C$ and subsequently stabilises at what they called the *lower buckling load, $P_L$*. They presented this as a useful empirical 'lower bound' for the collapse load of real shells with inevitable initial imperfections and finite disturbances. It was next shown by Friedrichs [1941] that just above $P_L$ there was a load, $P_M$, at which the grossly deformed but stabilized periodic path first had a total potential energy less than that of the trivial state. Tsien [1942] seized on $P_M$ as the *energy criterion load*, comparing it to Maxwell's rule in thermodynamics. He implied, illogically and without any evidence, that this would be a logical failure load due to dynamic disturbances. This explanation as it stands could not be true, since it would imply that near $P_M$ the shell would be constantly jumping in and out of its buckled state (unless the first inwards jump damaged the material of the shell).

A more complete survey of the history of shell buckling, looking particularly at features relevant to localization and the Maxwell load, is given in an Appendix.

More recently, following the outline of the static-dynamic analogy by Hunt, Bolt and Thompson [1989], it has been used with great success by Hunt and his co-workers in Bath and Bristol [Hunt & Lucena Neto,1993; Hunt, *et al*, 2000; Budd, *et al*, 2001; Hunt *et al*, 2003; Hunt, 2011], who looked in particular at the cylindrical shell. In this work the Maxwell load emerged naturally as a major pattern-formation feature in the large-deflection post-buckling response of shells and other 'long thin' structures. Notably, it was established that above the Maxwell load there exists a finite amplitude post-buckling path of unstable spatially localized states. For the cylindrical shell, this work is complicated by wave-number changes and chaotic solutions.

A particularly significant paper by Horak, Lord and Peletier [2006] on the mountain pass as an organizing centre in cylinder buckling demonstrated that the localised solutions offer a dramatically reduced energy barrier against lateral disturbances, be they static or dynamic. It is this important finding that we seek to clarify and extend in the present paper, drawing on analogies with the buckling of twisted rods.

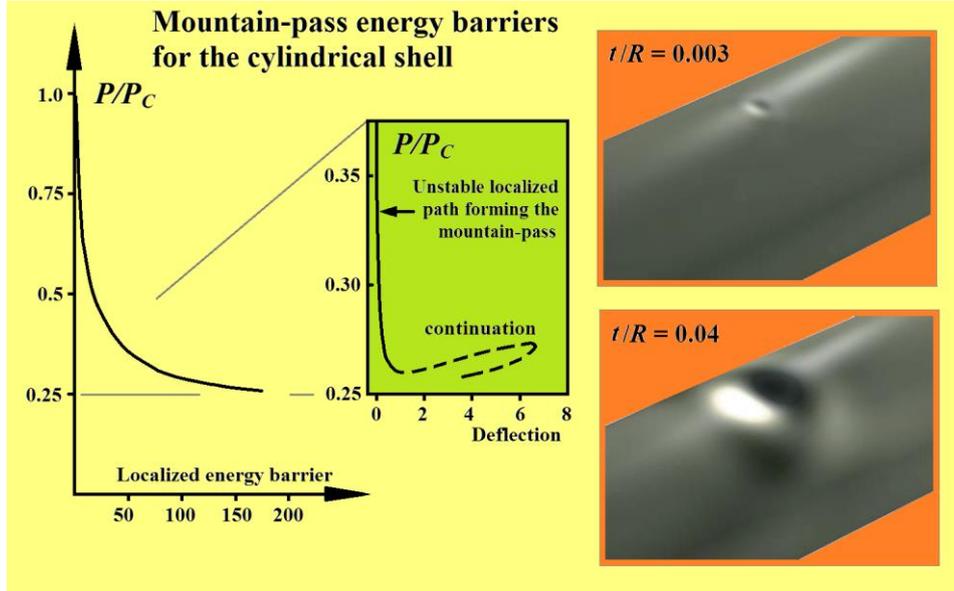

Fig. 1. Summary of the results of Horak, Lord & Peletier [2006] on the buckling of an axially loaded cylindrical shell.

## 2. Scope of the present paper

This paper draws on, and extends, the significant recent work at Bath and Bristol on the static-dynamic approach to the buckling and post-buckling of 'long thin' structures including struts on nonlinear elastic foundations and axially loaded cylindrical shells. This material, when combined with a stream of research at University College London (UCL) on the spatial buckling of stretched and twisted rods, allows new features of the 'energy criterion' Maxwell load to be observed and assessed. Using, in particular, an integrable problem of a rod constrained to lie on a cylinder [van der Heijden, 2001] fine details of the heteroclinic connection between the trivial unbuckled solution and the re-stablilized large-amplitude periodic state are examined and explained. It is this connection at the Maxwell load that destroys the falling spatially-localized post-buckling equilibrium path that descends from the classical linear buckling load.

This unstable localized path offers an order-of magnitude lower structural energy barrier (mountain pass) than the unstable periodic path, as elegantly demonstrated by Horak, Lord and Peletier [2006] for the cylindrical shell. So the practical significance of the Maxwell load is that it signals the onset of a super-sensitivity to static and dynamic lateral disturbances that we shall refer to as *shock sensitivity*. We demonstrate how this appears instantaneously at the Maxwell load when the high barrier of the periodic path is suddenly superseded by the very low barrier of the localised path, and show that its magnitude (which we write as $E^*$) can be calculated by an integration along the heteroclinic phase arc.

We should finally note at this point that throughout the paper we adopt the historical definition of the Maxwell load as that at which another (post-buckled) solution acquires lower energy than the trivial solution. Very often, the relevant solution which acquires a structural potential energy equal to that of the trivial fundamental path will be a re-stabilized periodic (helical for a rod) post-buckling path. In all the problems that we examine here this periodic path is unique. It then remains, for example, to show (as in §4) that it is at this Maxwell load that the localizing path falling from the classical critical bifurcation point is destroyed in a heteroclinic saddle connection.

## 3. Equivalent oscillators and localization

### 3.1 Strut on a nonlinear foundation

Drawing on the work of Potier-Ferry [1983] described in the Appendix, Hunt *et al*. [1989] employed the static-dynamic analogy in the study of structural localisation using in particular the beam (strut) on a nonlinear elastic foundation described by the ordinary differential equation (ODE)

$$EIy'''' + Py'' + ky - cy^m = 0 \quad \text{with } m = 2 \qquad (1)$$

Here $EI$ is the bending stiffness of the beam, $P$ is the axial compressive load at the ends, $k$ is the linear foundation force per unit length and $c$ is the magnitude of the foundation's nonlinearity. The lateral deflection of the beam is $y(x)$ where $x$ is the independent length along the beam. The index, $m$, is taken as 2 in the analysis, compared with the $m = 3$ used by Potier-Ferry. The above equation is technically 'non-integrable' so, like the cylindrical shell, it gives rise to spatial chaos: much fundamental work has gone into exploring its full behaviour [Buffoni, *et al*, 1996].

The total potential energy of this structure, incorporating the strain energy of bending, the foundation energy and the energy of the dead compressive load, $P$, can be written as

$$V_{ST} = \tfrac{1}{2} EI \int y''^2 \, dx + \int [\tfrac{1}{2} ky^2 - (1/3) cy^3] \, dx - \tfrac{1}{2} P \int y'^2 \, dx \qquad (2)$$

where the subscript 'ST' denotes *structure*. Assuming the beam to be very long (strictly infinite) linear theory gives the critical buckling load and buckling mode as

$$P_C = 2 \sqrt{(kEI)} \quad \text{and} \quad y = \cos \omega x \quad \text{where} \quad \omega^2 = \sqrt{(k/EI)}. \qquad (3)$$

The essence of the static-dynamic analogy is to regard the system as a *virtual* dynamical system by simply identifying $x$ as the time $t$. We use the word *virtual*, here, because no attempt is made to identify or devise a mechanical realisation of the dynamical system (in contrast to the twisted rod that has as its exact and real dynamical analogy, the spinning top).

The root structure of the ODE, linearized about the trivial state $y = 0$, illustrated in figure 2, is now seen to be what in non-linear dynamics would be called a Hamiltonian Hopf bifurcation, This is observed, for example, in the flutter of an undamped cantilevered pipe discharging a fluid.

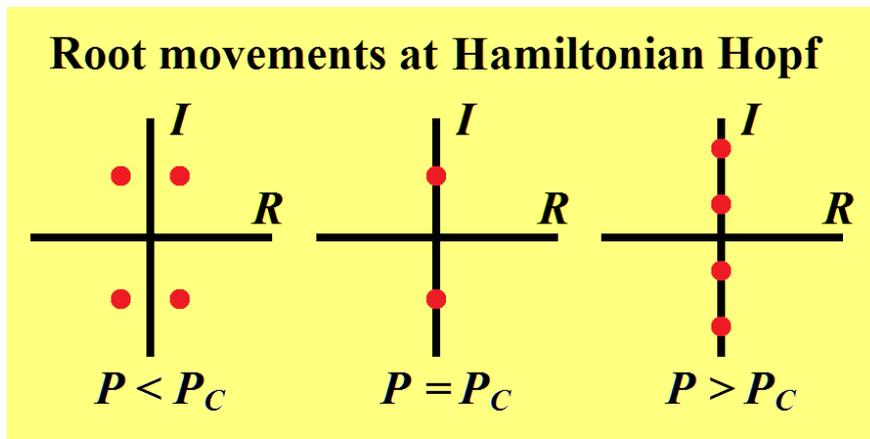

Fig. 2. The root movements typical of a Hamiltonian Hopf bifurcation.

A double-scale perturbation analysis, writing

$$y = y_1(x)\,\epsilon + y_2(x)\,\epsilon^2 + y_3(x)\,\epsilon^3 + \ldots, \tag{4}$$

gives the localizing solution of interest, valid close to the bifurcation point as, to successive orders in epsilon,

First perturbation equation: $\quad y_1 = A(X)\cos\omega x \tag{5}$

Second perturbation equation: $\quad y_2 = cA^2(9 + \cos 2\omega x)/18k \tag{6}$

Third perturbation equation: $\quad 2P_C\,A_{XX} - \omega^2 A + (19c^2/18k)A^3 = 0 \tag{7}$

Here the 'envelope' amplitude, $A(X)$, is a function of the 'slow' independent variable, $X := \epsilon x$, where the small parameter $\epsilon := \sqrt{(P_C - P)}$ is a measure of the distance from the critical bifurcation point; and we are writing $A_{XX} := \partial^2 A/\partial X^2$.

The third equation, thinking now in terms of the times $(t, \tau)$ rather than the distances $(x, X)$, is that of a one-degree-of freedom nonlinear oscillator of displacement $A(\tau)$, corresponding to a particle moving in a potential well of the form

$$V_{EO} = -\tfrac{1}{2}\omega^2 A^2 + \tfrac{1}{4}(19c^2/18k)A^4. \tag{8}$$

Here the subscript 'EO' denotes *equivalent oscillator*, and this potential energy should not be confused with the potential energy of the structure, $V_{ST}$.

This twin-well symmetric potential, $V_{EO}(A)$, has a maximum at $A = 0$ corresponding to the *stable* beam in its trivial un-buckled state, and minima at $A = \pm(\omega/c)\sqrt{(18k/19)}$ which correspond (with fixed $A$) to an *unstable* periodic post-buckled state emerging from the pitchfork bifurcation at $p := P/\sqrt{(kEI)} = 2$. Meanwhile the saddle-loop from the trivial solution defines the localized path as illustrated in figure 3. Note that the familiar rule, stable = minimum and unstable = maximum, is reversed when dealing with $V_{EO}$.

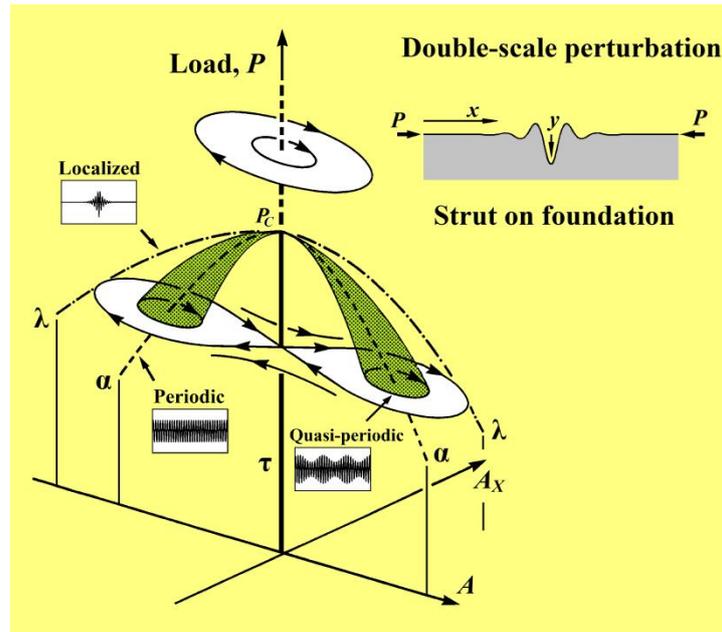

Fig. 3. Response of the equivalent oscillator in a double-scale perturbation analysis of a beam on a nonlinear elastic foundation

At the time of writing, Hunt, Bolt and Thompson [1989] did not evaluate the energy barriers associated with the two unstable post-buckling paths, but Giles Hunt has kindly made the required calculations for us in October 2013! The results are shown in figure 4, comparing the barriers of the displayed localized solutions with the barriers offered by the drawn (multi-wave) periodic forms. The solid circle shows the contribution of a single periodic wave at $p = 1.5$. We shall be discussing the features of this, and similar barrier-versus-load graphs, throughout the paper.

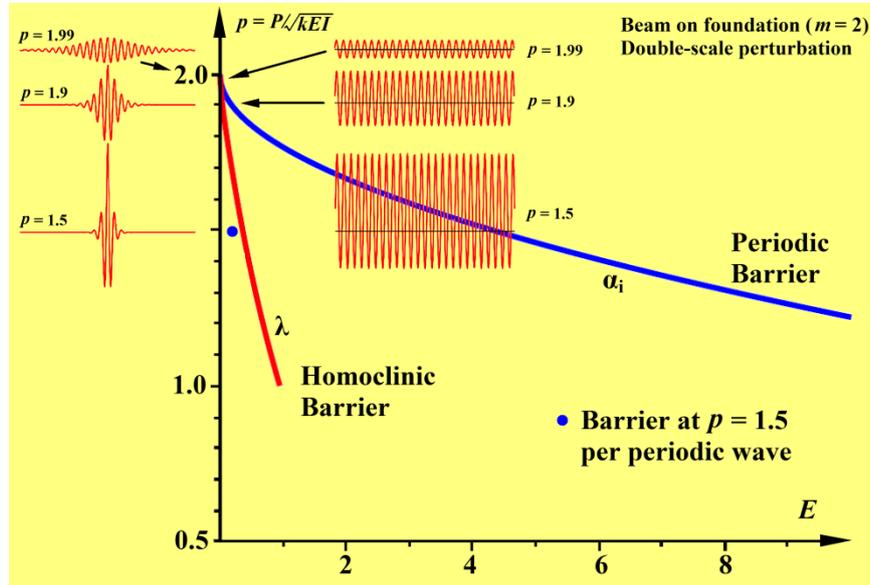

Fig. 4. Mountain pass energy barrier for the beam on a nonlinear elastic foundation (kindly calculated and supplied by Giles Hunt in 2013)

### 3.2. Stretched and twisted rod

A useful archetypal model for studying localization phenomena is the long (mathematically infinite) spatially-deforming twisted elastic rod subjected to an end tension and end twisting moment. This was studied theoretically and experimentally (Fig. 5) by Thompson and Champneys [1996].

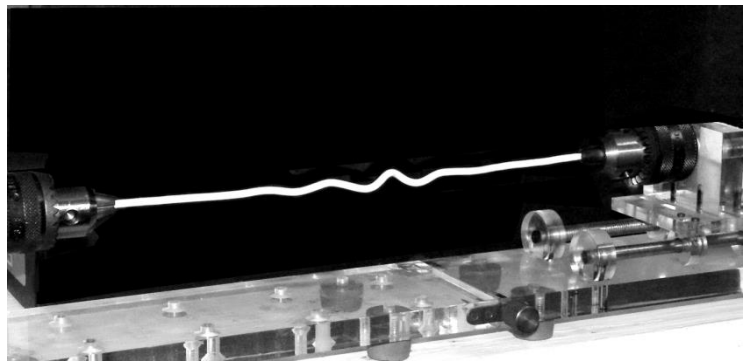

Fig. 5. A localized form in the post-buckling of a stretched and twisted elastic rod [Thompson & Champneys, 1996]. The silicon rubber rod is loaded rigidly by a controlled end displacement and a controlled end rotation. This rigidity allows the observation of what, under dead loading, would be an unstable state.

As first noted by Kirchhoff, the statics of the isotropic (or anisotropic) stretched and twisted rod is mathematically equivalent to the dynamics of a symmetric (or non-symmetric) spinning top if we identify the axial (arc-length) coordinate of the rod, $s$, with the time variable of the top, $t$. As with all dynamic analogies of elastic solids under conservative loading, the dynamical problem is itself conservative, with no damping.

In comparing tops and rods, the following notations and equivalences are used. The three principal moments of inertia of the top, ($I_1$, $I_2$, $I_3$), are equivalent to the two principal bending stiffnesses of the rod ($B_1$, $B_2$) and the torsional stiffness, $C$. The independent variables are the time, $t$, and arc length, $s$, while the dependent variables are the three Euler angles ($\theta$, $\psi$, $\phi$). For the rod, which is our main concern, these are the angles of a triad of vectors that we imagine spiralling along the spatially twisted centre line, with $\theta = 0$ corresponding to the loading axis. The two loads are an axial tension, $T$, and an end twisting moment, $M$. These can both be combined in the analysis to give the dimensionless loading parameter

$$m := M/\sqrt{(BT)} \qquad (9)$$

In the symmetric case [van der Heijden & Thompson, 2000] when the top has circular symmetry ($I_1 = I_2 = I$) and the rod (having, for example, a circular cross-section) has equal principal bending stiffnesses ($B_1 = B_2 = B$), there are three conserved quantities. These are, for the top (followed by the rod), the Hamiltonian (ditto), the angular momentum (torque) about a vertical axis (the loading axis), and the angular momentum (torque) about the body symmetry axis. These three conserved quantities make the symmetric problems integrable, so no (dynamic or static) chaos can be encountered. Meanwhile in the general problem ($I_1 \neq I_2$, $B_1 \neq B_2$) the body-axis terms are no longer conserved, the integrability is lost, and chaos can be found. A rigorous proof of chaos in the non-symmetric top, and by extension the asymmetric rod, has been presented recently [van der Heijden & Yagasaki, 2013]. This improves on earlier treatments, such as [Holmes & Marsden, 1983], which did not take sufficient account of a singularity of the equations typical of Euler-angle formulations.

The analysis of the symmetric rod proceeds as follows. In terms of the Euler angles $\theta(s)$, $\phi(s)$ and $\psi(s)$, the total potential energy of the structural system, being the sum of the rod's strain energy, and the work contributions from the end tension $T$ and end twisting moment $M$, is given by:

$$V_{ST} = \int (\tfrac{1}{2}B\kappa^2 + \tfrac{1}{2}C\tau^2)\,ds + T\int (1 - \cos\theta)\,ds - M\int (\phi' + \psi')\,ds \qquad (10)$$

where a prime denotes differentiation with respect to $s$, and the curvature, $\kappa$, and kinematic twist, $\tau$, are given by

$$\kappa^2 = \theta'^{\,2} + \psi'^{\,2}\sin^2\theta \qquad (11)$$

$$\tau = \phi' + \psi'\cos\theta. \qquad (12)$$

Introducing the Lagrangian (total potential energy density), $L$ we can write

$$V_{ST} = \int L(\theta, \theta', \psi, \psi', \phi, \phi')\,ds \qquad (13)$$

The Euler-Lagrange equations for $\psi$ and $\phi$ are easily solved (the latter giving $\tau = M/C = $ constant), and back substitution gives

$$V_{ST} = \int L(\theta, \theta')\,ds \qquad (14)$$

in terms of the reduced Lagrangian

$$L(\theta, \theta') = \tfrac{1}{2}B\theta'^{\,2} - (M^2/2B)(1 - \cos\theta)/(1 + \cos\theta) + T(1 - \cos\theta) - M^2/2C \qquad (15)$$

Taking a dynamical systems view, we now think of this as an equivalent oscillator, with

$$L = K - V_{EO} \tag{16}$$

where the kinetic energy of the equivalent oscillator, $K$, is

$$K = \tfrac{1}{2}B\theta'^{2} \tag{17}$$

and the potential energy of the equivalent oscillator, $V_{EO}$, is

$$V_{EO} = (M^{2}/2B)(1-\cos\theta)/(1+\cos\theta) - T(1-\cos\theta) + M^{2}/2C \tag{18}$$

Now, $K + V_{EO}$ = const. gives us contours corresponding to the phase portrait of the 'oscillator'.

So the statics of the rod is governed by the equivalent one-degree-of-freedom nonlinear oscillator in the $(\theta, \theta')$ phase space, behaving as a particle on a total potential energy surface given by $V_{EO}(\theta)$. This is analogous to the perturbation solution of the beam on an elastic foundation, but is here the complete solution valid at any distance from the governing bifurcation. From now on we can focus on the previously defined dimensionless control parameter, $m := M/\sqrt{(BT)}$, which is the only system parameter if we ignore the 'arbitrary' constant $M^{2}/2C$ (which would introduce the ratio of $B$ to $C$). The energy function, $V_{EO}(\theta)$, varies with $m$ in such a way that a sub-critical pitchfork bifurcation occurs at $m = 2$. For $m < 2$ the function has two non-trivial minima, as shown in figure 6, while for $m > 2$ there is just a single minimum at $\theta = 0$. Two phase portraits of the oscillator are sketched in white.

Thinking first about the top, the flow along these phase trajectories represents its dynamical motions with the time, $t$. Here we must think of our composite control parameter, $m$, as a measure of its speed of spin. For $m > 2$ the upright top is stable, and perturbations of its trivial state are governed by the closed obits around the origin $(\theta, \theta') = (0, 0)$. Meanwhile for $m < 2$ the upright top is unstable, and small perturbations will send it off to finite deflections. Keeping close to the drawn homoclinic loop, it will nevertheless return periodically, on a long time scale, to the near-upright state.

For the rod, however, a particular phase trajectory represents the variation of its static deformation as we progress along its centre-line; and the fixed point at the origin corresponds to the straight (but twisted) rod. Unlike the top, the trivial undeflected configuration of the rod is stable for loads corresponding to $m < 2$, and unstable for $m > 2$.

The linear eigenvalue buckling load of the rod is given by $m = 2$ from which is generated a sub-critical pitchfork bifurcation. For $m > 2$ the only fixed point is the trivial fixed point. Meanwhile at values of $m < 2$ for the infinitely long rod (with no imposed end conditions) we have in addition to the stable fundamental solution a sub-critical helical path corresponding to the non-trivial fixed points of the oscillator, modulated helices corresponding to the closed orbits, and finally a localised solution corresponding to the homoclinic saddle loops which start and end arbitrarily close to the trivial state. Both helical and localized paths bifurcate from the trivial un-buckled solution at the critical buckling load, $m_C = 2$. They sustain a falling load (exhibiting no $P_L$ or $P_M$), and are unstable under dead loading.

The thumb-nails in the grey frame in figure 6 illustrate the static localized shapes of the rod, and as $m$ is decreased steadily from $m_C = 2$, the localization becomes increasingly pinched. Notice that if we were to impose the condition that the rod line up with its original axis at plus or minus infinity, only the localised forms would be preserved.

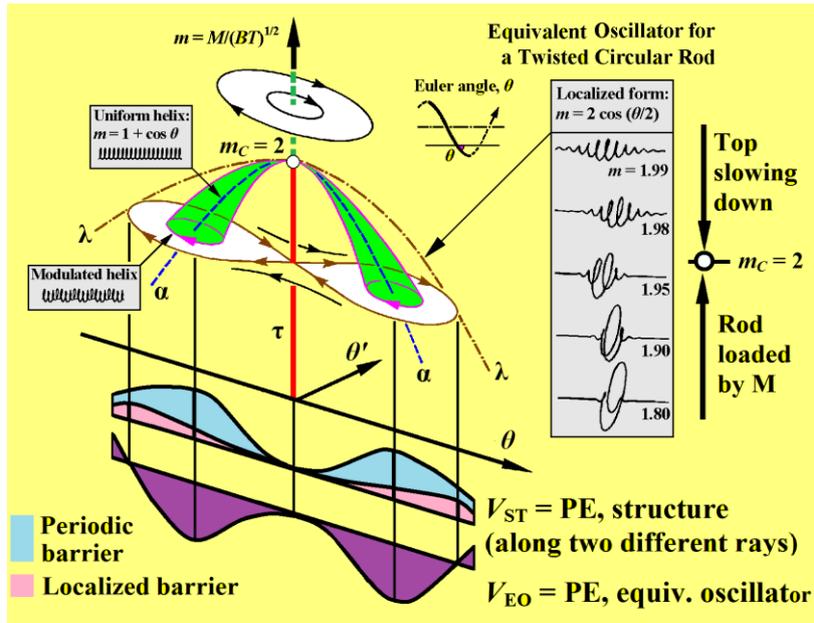

Fig. 6. Phase portraits and potential energy plots for the stretched and twisted rod

At the bottom of this figure we sketch schematically the total potential energy of the equivalent oscillator, $V_{EO}(\theta)$, for the sub-critical phase portrait at $m < 2$. This should not be confused with the total potential energy, $V_{ST}$, of the structure (the rod and its dead loads) which is a function of the shape of the rod, not just $\theta$. We represent this very schematically by two curves to show how it might vary along notional periodic or localizing rays from the trivial solution in the deformation space (see figure 9 for a more realistic 3D representation).

The question now arises as to which unstable solution offers the lower energy barrier resisting disturbances of the metastable trivial solution. To examine this we have extended the previously published solution [van der Heijden & Thompson, 2000] to obtain the energy barrier diagram shown in figure 7. This shows that, as with the beam on an elastic foundation, the energy barrier offered by the localized ($\lambda$) form is roughly equal to the energy of one periodic helical wave.

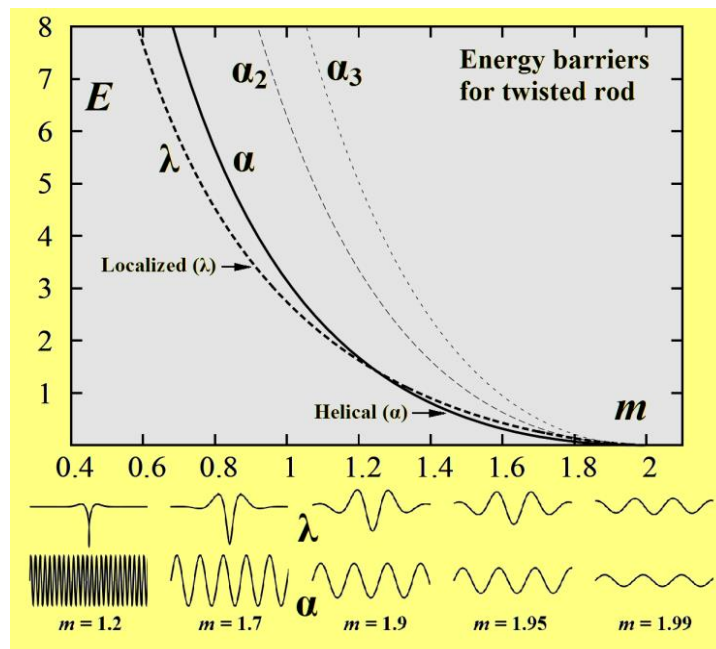

Fig. 7. Energy barriers and thumb-nails of the deflected forms for the stretched and twisted rod.

## 3.3. The twisted rod constrained to a cylinder

We have seen the significance of a localized post-buckling path in offering a reduced energy barrier against lateral disturbances, and we next need to examine how this path is created or terminated (under increasing or decreasing load) at the Maxwell condition.

A particularly useful integrable illustration is offered by the behaviour of a long stretched and twisted *isotropic* rod constrained to lie on a cylinder (akin to a drill-string, when the constraining pressure is inwards) which was studied by van der Heijden [2001]. This rod on a cylinder has the same composite loading parameter as the earlier free rod, $m := M/\sqrt{(BT)}$, where $M$ is the applied twisting moment, $T$ is the applied tension, and $B$ is the bending stiffness about any axis. Its exact solution is, once again, governed by an equivalent one-degree-of-freedom oscillator. It has an infinite 'buckling' load, $m_C$, with an unstable falling helical post-buckling path, together with a falling localizing solution. In a certain parameter range, the helical path restabilizes at a lower buckling load, $m_L$, with a Maxwell load (defined by the equality of trivial and stable-periodic energies) at $m_M > m_L$.

So for the free twisted rod [van der Heijden & Thompson, 2000] we have a finite $m_C$ but no $m_L$ or $m_M$ while for the rod on the cylinder [van der Heijden, 2001] we have a finite $m_L$ and $m_M$, but $m_C$ is infinite. Drawing on these, the schematic *composite* diagram of figure 8 illustrates a general case (not necessarily a twisted rod), in which all three critical loads are finite. Notice that when talking in general terms about a 'long thin structure', which might be a rod or a shell, etc, we use the expression 'periodic' rather than 'helical' (the latter is specific to rods), and we call the load $P$.

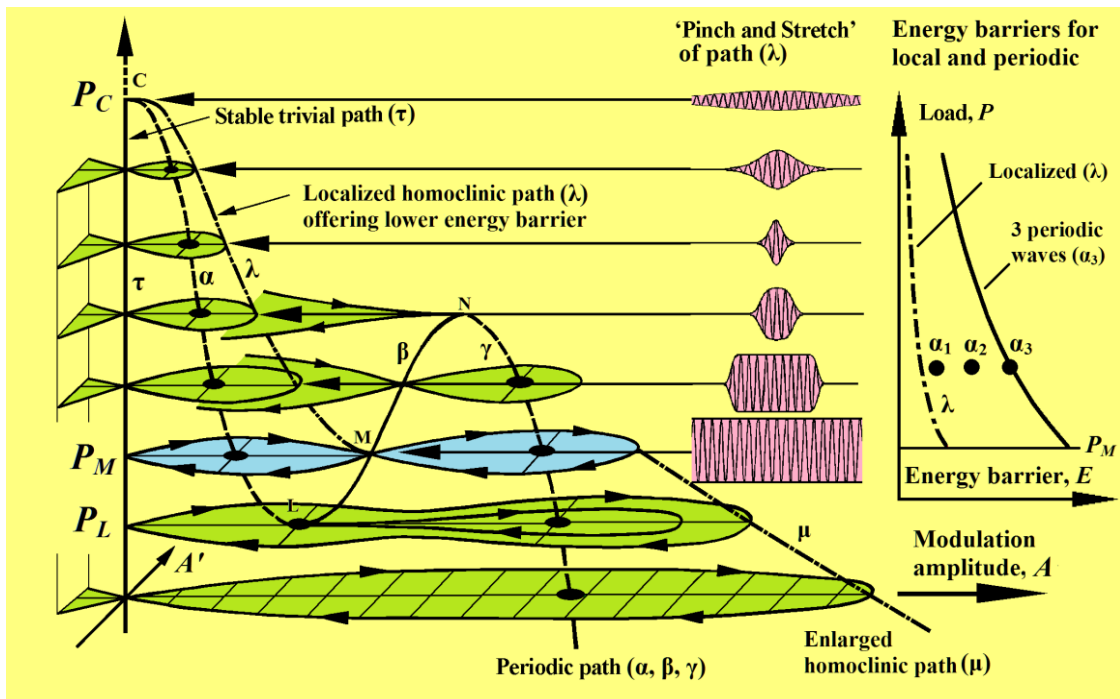

Fig. 8. A schematic diagram illustrating the buckling of a shell-like structure. The 2D phase portraits of the equivalent oscillator are sketched in the space of the modulation amplitude, $A$ and its time derivative $A'$ at various fixed values of the load, $P$. Here M identifies the collision between the two post-buckling paths at the Maxwell load. Diagrams on the centre-right show the current shape of the localizing mode. The load versus energy barrier plots on the far right relate to the rod on a cylinder.

Notice that to allow a concise discussion, and easy comparison between different figures, we have introduced the following path nomenclature in terms of the bifurcation point, C, the fold at the lower buckling load, L, any subsequent fold of the periodic path, N, and the point of collision between the homoclinic path and the periodic path, M:

Path α: unstable periodic path from C to L (for the energy barrier, α is for 1 period, $α_2$ is for 2, etc)
Path β: stable periodic path after L and up to any subsequent folding at N
Path γ: unstable periodic path after N
Path λ: the localizing homoclinic path from C to M
Path μ: the enlarged homoclinic path after the 'explosion'
Path τ: the stable trivial fundamental path from the unloaded state to C

The equivalent oscillator varies with the controlled load parameter, *P*, and sketches are made of the two-dimensional dynamical phase spaces spanned by the modulation amplitude, *A*, and its space/time derivative *A′* at different values of *P*. As before, stable (unstable) states of the structure appear as unstable (stable) states of the equivalent dynamical oscillator. In the following, stability statements refer to the rod or shell. The non-trivial fixed points of the oscillator, shown as small black ellipses, define the periodic post-buckling path which is unstable (α) down to $P_L$ rises (β) to a maximum and then falls again (γ). Meanwhile the saddle connections from the trivial solution (τ) on the *P* axis define the localized homoclinic path (λ) which falls until it is destroyed on colliding with the stable rising branch (β) of the periodic path at $P_M$. Sketched on the centre-right are the notional deformed shapes of the structure as the load decreases along the homoclinic path (λ), forming what we have called the 'pinch and stretch' effect [van der Heijden, *et al*, 2002].

In this paper [van der Heijden, *et al*, 2002], the generalization to an *anisotropic* rod on a cylinder is non-integrable and generates spatially chaotic deformations; but this paper is nevertheless a useful reference for the formulation of the *isotropic* rod. In particular, Fig. 14 of [van der Heijden, *et al*, 2002] shows for the isotropic rod the stabilization of the periodic post-buckling path at $P_L$, together with the localized snaking generated by the anisotropy.

## 4. Path collisions at the Maxwell load

The work on the symmetric rod on a cylinder [van der Heijden, 2001] is currently being extended [van der Heijden & Thompson, 2014] to look in detail at conditions near the Maxwell load. This current work on the integrable system gives the clearest view yet of the detailed collisions and energy properties at the Maxwell load, and we give here some of our preliminary findings. Firstly, we establish analytically that the path collision at M does indeed correspond to the Maxwell load as defined historically by $V_{ST}(τ) = V_{ST}(β)$. Secondly, we observe the details displayed in figure 9.

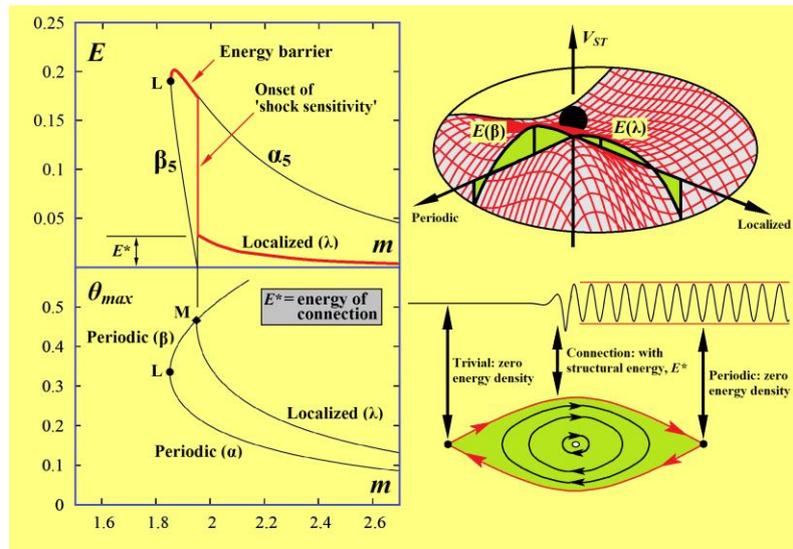

Fig. 9. This figure summaries the findings of current work to be published [van der Heijden & Thompson, 2014]. The composite graphs on the left hand side show the evaluated energy hump *E(m)* and the equilibrium paths as $θ_{max}(m)$. Top right is a 3D sketch of the structural energy $V_{ST}$, showing the periodic and localized energy barriers. Finally, the lower right diagram illustrates how the localized energy barrier at the Maxwell load can be evaluated via an arc of the heteroclinic connection.

The lower left diagram shows on a plot of $\theta_{max}$ against control parameter, $m$, the equilibrium paths, β and λ colliding at the Maxwell point, M. Meanwhile the upper left diagram shows the evaluated structural energy barrier, $E(m)$, corresponding to $V_{ST}$ measured from the datum of the trivial path. Here $E(\lambda)$ for the localized path is compared with the energy barrier for five periods of the periodic path, $E(\alpha_5)$. Notice that if we had plotted the energy for a single period, $E(\alpha)$, it would lie very close to $E(\lambda)$, following what we have observed for the beam on a nonlinear foundation and for the unconstrained twisted rod. We have omitted it here, firstly for clarity, and secondly because a single period (being on its own kinematically inadmissible in the sense that it must be part of a continuous wave) does not offer a valid escape route for comparison purposes. The choice of 5 periods is of course rather arbitrary, but any structure deemed to be 'long' must accommodate 5 periods as an absolute minimum. As always, in dynamic-analogy discussions, we are ignoring the supposedly remote boundary conditions. The red line emphasises the sudden loss of safety barrier at $m_M$, defining the onset of *shock sensitivity*.

The upper right hand diagram simply illustrates the form of the multi-dimensional $V_{ST}$. Meanwhile the lower right diagram illustrates an interesting feature of the $E(m)$ graphs. At the heteroclinic connection, in which path λ collides with path β, the localized and periodic solutions do not have the same $V_{ST}$ or the same $E$ (recall, here, that $E$ differs from $V_{ST}$ only by being measured from the datum of the trivial solution). Thus while path β has $E$ equal to zero, in line with the Maxwell definition, path λ ends at $m_M$ with the structural energy barrier $E^*$. The latter can be obtained analytical by integration around the heteroclinic half-path in phase space.

## 5. Practical significance of the Maxwell load

We have seen in a number of examples that the path of localized small-amplitude states (existing only above the Maxwell load) presents a total potential energy barrier roughly equal to the energy of a single periodic wave. The difference is, of course, that the localised form is a kinematically admissible displacement, while the periodic form must embrace a large number of waves. This means that the energy barrier of the localized form is an *order of magnitude* lower than the energy barrier of the periodic solution. Indeed, as the length of the structure goes to infinity, so does the periodic energy barrier.

This then is the practical significance of the Maxwell load, its creation of a low-amplitude localized path which gives the first realistic possibility of an escape from the trivial state due to dynamic disturbances or static side loads. This onset of a *super-sensitivity* to lateral disturbances can be summarized as:

> ***The Maxwell load signals the onset of 'shock sensitivity'.***

This result can have significance in many areas of applicable mechanics and in the emerging field of Turing pattern formation [Dawes, 2010].

## 6. Discussion and conclusions

We have shown the true significance of the Maxwell load in shell-like buckling in terms of 'shock-sensitivity'. So how does this impinge on the much studied scatter of experimental buckling loads for the axially compressed cylindrical shell, shown in figure 10?

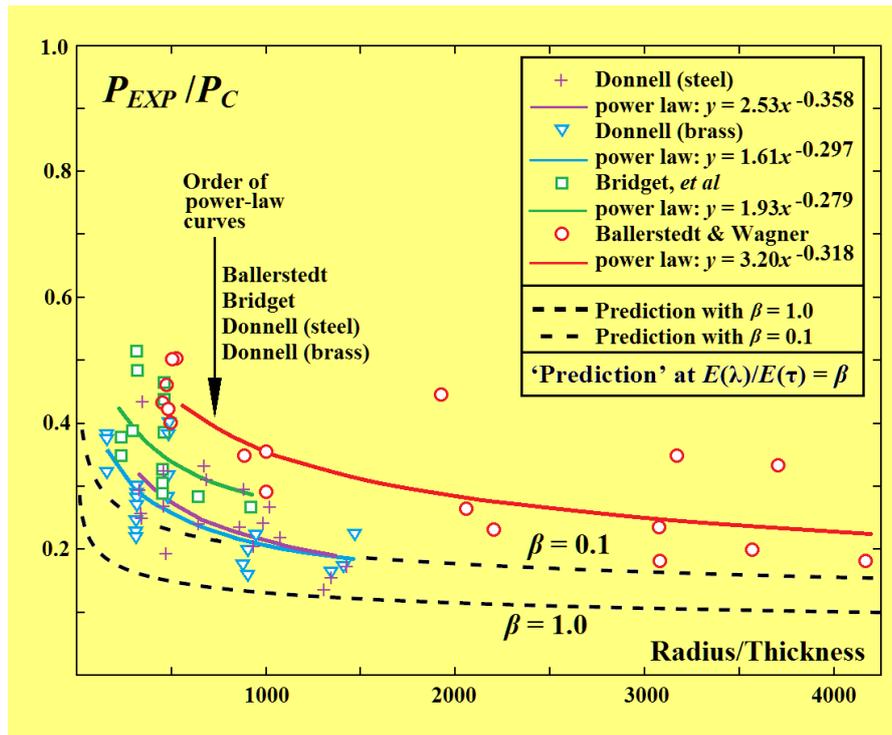

Fig. 10. Historical experimental buckling loads for the axially compressed cylindrical shell adapted from figure 5.3 of Horak, Lord and Peletier [2006]. The theoretical 'lower bound' estimates based on the localized mountain pass analysis of these authors are shown as dotted lines for two values of $\beta$, namely the ratio of the energy barrier to the stored energy in a cylindrical section of length $2\pi R$, where $R$ is the radius of the shell. Their power law indexes, such as their 0.358, compare with the 0.5 that Zhu, et al [2002] derive from a larger set of data points.

It is unlikely that such lateral shocks were present in careful laboratory tests, so the answer would seem to lie in the enhanced sensitivity to *localized imperfections* in the shells that would surely accompany the demonstrated super-sensitivity to lateral disturbance. This is something that should certainly be addressed by our mathematical colleagues, embracing, perhaps, the probabilistic viewpoint of Elishakoff [1983]. Meanwhile, Calladine and his co-workers [Zhu, et al, 2002] addressed the statistical distribution of the experimental failure loads. They pointed to the fact that classical buckling loads (and other derived loads, such as $P_M$ and $P_L$) are proportional to the thickness of a cylindrical shell, while the experimental results scale as thickness to the power 1.5. They present no new theory to remedy this anomaly, but argue the case for the 1.5 power from a variety of theoretical and experimental observations. Notice that for the form of presentation used in figure 10, Calladine's 1.5 would correspond to an index of 0.5 (compared with the 0.358 shown for one of the data sets). One wonders whether the 1.5 could perhaps be explained by the scaling of manufacturing errors with the thickness to radius ratio.

 The drawn 'predictions' used by Horak *et al* [2006] to offer a notional lower bound to the experimental results are drawn for two values of $\beta = 0.1$ and $\beta = 1.0$. These are obtained by using the hypothesis that failure can be expected once the energy barrier has dropped to $\beta$ times the stored energy in the trivial configuration of a cylindrical section of length $2\pi R$ where $R$ is the radius of the shell. Note that this approach gets away from the Calladine objection that derived loads are proportional to the shell thickness.

 The new feature of the Maxwell load described in this paper, involving an instantaneous and order-of- magnitude reduction in the escape barrier against lateral disturbances, has applications that lie well beyond the traditional field of thermodynamics. It is significant for practical thin shells in *operational* environments, and underlies many important phenomena of localization in solids and fluids including pattern formation by the Turing instability in chemical and biological kinetics [Dawes, 2010]. Another new application may arise from a most significant recent study by Cao and Hutchinson [2012] of the surface wrinkling of a compressed half-space. They show that the wrinkling is both highly unstable and extremely imperfection-sensitive with close parallels to the behaviour of the cylindrical shell under axial compression. Indeed the mathematics underlying the nonlinear wrinkling behaviour is very similar to that discovered by Koiter for cylindrical shell buckling. Their analysis suggests that a near-perfectly flat surface would collapse dynamically into a creased state before appreciable wrinkling can be observed. This is entirely analogous to the cylindrical shell snapping dynamically into the collapsed state without revealing its short wavelength mode.

**Appendix: History of the Maxwell load in shell buckling**

The buckling of thin elastic shells is characterised by a very unstable, sub-critical post-buckling path of periodic equilibrium states, which falls rapidly from the classical linear buckling load, $P_C$, and subsequently stabilises at what was historically called the *lower buckling load, $P_L$*. Experiments, motivated in particular by load-carrying aerospace shells, had established that uniformly compressed spherical shells and axially loaded cylindrical shells failed at loads which could be as low as ¼ $P_C$, and von Karman & Tsien [(1939, 1941] suggested that $P_L$ could be adopted as a useful empirical 'lower bound' for the collapse load of real shells with inevitable initial imperfections and finite disturbances. For the spherical shell, $P_L/P_C$ was subsequently shown experimentally to be about 0.25 for a polythene sphere which could be loaded and unloaded in the elastic range; this agreed accurately with the theoretical post-buckling solution [Thompson, 1962]. A fairly similar ratio was found by Yamaki [1984], in his experimental post-buckling studies of the cylindrical shell where changes in circumferential wave number complicate the story.

It was next observed by Friedrichs [1941] that just above $P_L$ there was a load, $P_M$, at which the grossly deformed but stabilized periodic path first had a total potential energy less than that of the trivial state. He did not suggest this as a buckling criterion, but rather as a conveniently determined approximation to $P_L$.

Tsien [1942], on the other hand, seized on $P_M$ as the *energy criterion load*. Comparing it to Maxwell's rule that a thermodynamic system is likely to be found in its lowest energy state, he implied (incorrectly) that this would be a logical failure load due to dynamic disturbances. Given that the experiments were performed in careful laboratory tests where disturbances were minimal, this explanation could not be true, since it would imply that near $P_M$ the shell would be constantly jumping in and out of its buckled state.

Very shortly after these papers, Koiter [1945] demonstrated theoretically, in his classic dissertation, that thin shells had a severe sensitivity to initial geometrical imperfections. Written in Dutch, it took many years before this work became well-known, but when it did it seemed to many that it had solved the essence of the shell buckling problem; though some anomalies persist to the present day [Zhu, *et al*, 2002]. Next, Hoff *et al* [1966] made a fairly intensive theoretical study of the post-buckling of the axially compressed cylindrical shell, using a total potential energy expression consisting of about 1100 terms which was minimized with respect to 16 unknowns. They argued that the commonly-used von Karman-Donnell formulations were only valid for very low thickness-to-radius ratios, where the analysis found $P_L/P_C$ tending to zero. Progress in these formative years, from von Karman & Tsien to Hoff, was authoritatively assessed in *Applied Mechanics Reviews* by Hutchinson & Koiter [1970] and by Thompson & Hunt [1973].

A hint of other things to come was seen at the IUTAM Symposium at UCL in 1982 [Thompson & Hunt, 1983]. Here Tvergaard & Needleman [1983] gave an overview of localisation phenomena in a variety of structures including effects of plasticity and dynamics. Two other significant precursors studied localization of an infinitely long beam (or strut) on a nonlinear elastic foundation. Chater *et al* [1983] looked at dynamical buckle propagation where the Maxwell load has a key role to play. Potier-Ferry [1983] looked at localization phenomena in a compressed beam (strut) attached to a softening nonlinear elastic foundation with a cubic foundation force. Using a double-scale perturbation technique that had originated in fluid dynamics, he showed that the slow-scale *amplitude equation* was that of an equivalent oscillator. A homoclinic saddle-loop [Guckenheimer & Holmes, 1983; Thompson & Stewart, 1986] of this oscillator then described a sub-critical localizing path emerging from the classical bifurcation point at $P_C$. This was a seminal contribution, on which much progress has since been built, especially by Giles Hunt and his colleagues who have stimulated a renewed interest in the Maxwell load within the context of structural localisation studies.

Finally, we must say a few words about the different approaches of Isaac Elishakoff and Jim Croll to shell buckling. Elishakoff drew attention at the afore-mentioned IUTAM meeting to the need for a probabilistic approach to the buckling loads of imperfection-sensitive structures , a subject that he has pioneered over many years. He continues to argue with conviction [Elishakoff, 2012] that from an engineering standpoint the final target for the analyst must be to study imperfect shells (a challenge for many colleagues) and to express the results in probabilistic terms. Meanwhile Croll has developed an important analytical technique, see for example [Yamada & Croll, 1999; Croll & Batista, 1981], using his reduced-stiffness method which, by deleting certain energy terms in the formulation, aims to derive a lower bound to the experimental failure loads.

There is, of course, a massive literature on the analysis of imperfection-sensitivity in shell buckling that is beyond the scope of this brief article. The interested reader might like to take a look at the website of David Bushnell [2012].